\documentclass[11pt,mathsrc]{amsart}
\usepackage{amssymb,latexsym,epsfig}
\usepackage[mathscr]{euscript}

\makeatletter \@addtoreset{equation}{section} \makeatother

\newtheorem{Lemma}[equation]{Lemma}
\newtheorem{Theorem}[equation]{Theorem}
\newtheorem{Prop}[equation]{Proposition}

\newenvironment{Remark}{\noindent\textbf{Remark}}{}
\newenvironment{Proof}{\noindent\emph{Proof\ }}{\hfill$\square$\\}

\newcommand\QQ{\mathbb Q}
\newcommand\C{\mathbb C}

\newcommand\Z{\mathbb Z}

\newcommand\B{\,\overline{\!B}}
\newcommand\D{\,\overline{\!D}}
\newcommand\Zz{\,\overline{\!Z}}
\newcommand\Xx{\,\overline{\!X}}

\newcommand\Pee{\mathbb P^{\,1}}
\newcommand\OO{\mathscr O}
\newcommand\I{\mathcal I}

\newcommand\XX{\widehat X}

\title{Nodes and the Hodge conjecture}
\author{R. P. Thomas}
\date{}

\begin{document}
\maketitle \vspace{-5mm}
\begin{abstract} \noindent
The Hodge conjecture is shown to be equivalent to a question about the
homology of very ample divisors with ordinary double point singularities.
The infinitesimal version of the result is also discussed.
\end{abstract}

\section{Introduction}

Throughout this paper ``the Hodge conjecture" will mean the
statement that any rational class $A\in H^{2p}(X;\QQ)$ of pure
Hodge $(p,p)$ type on any smooth complex projective algebraic
variety $X$ is realised by a rational combination of
codimension-$p$ algebraic cycles in $X$. We show that this is
equivalent to a statement about the homology of nodal hyperplane
sections of even dimensional varieties. (By ``nodal" we mean a
hypersurface whose only singularities are analytically equivalent
to ordinary double points (ODPs). Such hypersurfaces are well
known to contain more of the middle dimensional homology of the
ambient variety than smooth hypersurfaces; think of reducible nodal curves
in surfaces, for instance, or in higher dimensions see for example
\cite{Cl}, \cite{Sch}.)

\begin{Theorem} \label{1}
The Hodge conjecture is true if and only if the following question
can be answered affirmatively for all even dimensional smooth complex
projective algebraic varieties $(X^{2n},\OO_X(1))$ and any class $A\in H^{n,n}(X;\C)
\cap H^{2n}(X;\QQ)$.

Is there a nodal hypersurface $D\subset X$ in $|\OO_X(N)|$ for some $N$,
such that $\mathrm{PD\,}[A]$ is in the image of
the pushforward map $H_{2n}(D;\QQ)\to H_{2n}(X;\QQ)$ ?
\end{Theorem}

(This can also be reformulated in terms of \emph{co}homology: Let $\D\stackrel
{\iota\,}{\hookrightarrow}\Xx\stackrel{\pi\,}{\to}X$ be the blow-up of
$D$ in its ODPs, inside the blow-up of $X$ at those points. Then the question
is whether $A$ is in the image of $\pi_*\iota_*:\,H^{2n-2}(\D)\to H^{2n}(X)$.)

This million dollar question appears, at first sight, to be much simpler than the
Hodge conjecture, asking only that we find a certain \emph{homology}
class, rather than an algebraic cycle (or even a pure $(p,p)$ homology
class);
i.e. that the cycle PD\,$[A]$ can be squeezed into a nodal hypersurface $D$
rather than fully into an algebraic $n$-cycle.
But it seems that finding nodal hypersurfaces may be nearly as difficult
as finding algebraic cycles, as we show infinitesimally in the last section
(though Theorem \ref{1} does at least bypass having
to find an algebraic cycle representing the above homology class in $D$,
which is itself a hard unsolved problem \cite{Sch}). Since the sections
of $\OO_X(N)$ with ODPs sit inside a vector space, and taken over all $N$
form a graded semigroup sitting inside the graded ring of $X$, perhaps
the more amenable tools of algebra might be
brought to bear on this formulation of the Hodge conjecture.

More likely, analysis and symplectic geometry might also be relevant. Using
Donaldson's work on Lefschetz pencils it is much better understood in
symplectic geometry how to contract vanishing cycles with a homology relation
between them to produce symplectic
nodal hypersurfaces containing extra homology classes; one could then try
to make such a hypersurface holomorphic by minimising its volume.
This should be a more manageable analytical problem then trying to minimise
the volume of a given middle dimensional homology class, as the cycle is
both symplectic and of low codimension. As
Herb Clemens pointed out to me, perhaps the interest in the theorem
is the converse, that the Hodge conjecture implies the existence of many
hyperplanes with many ODPs. Though it seems unlikely, known bounds on numbers
of such ODPs might then give a counterexample to the Hodge conjecture.
Similarly this might conceivably provide a symplectic route to disproving the
Hodge conjecture, by bounding numbers of Lagrangian spheres in smooth hyperplanes
of high degree; bounds like this exist in many low degree examples, whereas
the Hodge conjecture would imply there exist many such vanishing cycles
in high degree -- the vanishing cycles of a smoothing of the nodal hypersurface.
\\

\noindent \textbf{Acknowledgements.} Many thanks to Chad Schoen, Daniel
Huybrechts, Ivan Smith, Mark Gross and Miles Reid for useful conversations,
and to the Newton Institute, Cambridge, for its excellent research environment.
The author is supported by a Royal Society university research fellowship.

\section{Preliminaries}

Throughout we will denote by $\omega\in H^{1,1}(X)$ the first Chern class of the
polarisation $\OO_X(1)$, and by $H$ a divisor in its linear system, so
that $[H]=\omega$. First we will need to reduce the Hodge conjecture to the
following (which is presumably standard).

\begin{Prop} \label{2}
The Hodge conjecture is equivalent to the following statement.
Fix any even dimensional smooth complex projective algebraic variety $X^{2n}$,
and any class $A\in H^{n,n}(X;\C)\cap H^{2n}(X;\QQ)$.
Then there exist integers $N_1\ne0,\ N_2$ and an effective algebraic cycle $Z\subset
X$ whose fundamental class $[Z]$ equals $N_1A+N_2\omega^n$.
\end{Prop}

\begin{Proof}
If the Hodge conjecture is true, then any such $A$ is a rational
linear combination of cycles. That is, for some integers $N_i$,
$N_1A=N_3[Z_3]-N_2[Z_2]$, for effective algebraic $n$-cycles
$Z_i$.

$Z_2$ lies in an $n$-dimensional intersection of $n$ hypersurfaces
$H_i\in|n_iH|$ for sufficiently large $n_i$. (This is easily shown
inductively: any $n$-cycle in a (singular) $k$-dimensional variety
($k>n$) lies in a $(k-1)$-dimensional hypersurface of sufficiently
high degree.) That is, the intersection of hypersurfaces is
$Z_2\cup Z$ for some other effective algebraic $k$-cycle $Z$. Thus
$N_1A+N_2\big(\prod_1^nn_i\big)\omega^n=N_3[Z_3]+N_2[Z]$ is
effective, as required.

Conversely, we will show that if the given statement is true then
for any $k,\,d$ and $X^d$, and any rational class $A\in
H^{k,k}(X^d)$, there exist $N_1,\ N_2$ such that
$N_1A+N_2\omega^k$ is represented by an effective algebraic cycle.
This clearly implies the Hodge conjecture. \\

\noindent $\bullet$ For $k>d/2$, the hard Lefschetz theorem \cite{GH} gives
a rational class $A'\in H^{d-k,d-k}(X)$ such that $A'.\,\omega^{2k-d}=A$.
If we can find an effective algebraic cycle $Z'$ such that $[Z']=N_1A'+
N_2\omega^{d-k}$, then its intersection $Z$ with $(2k-d)$ generic hyperplanes
will be an effective cycle satisfying $[Z]=N_1A+N_2\omega^k$. Thus it
is sufficient to prove the statement for $A'$, so we may replace $A$ by $A'$
and $k$ by $d-k$. \\

So we can assume that $k\le d/2$, and proceed by induction on $d$ (for
fixed $k$); the base case $d=2k$ is the statement in the Theorem (with
$k=n$ and $d=2n$) which we are assuming true. \\

\noindent $\bullet$ For $k<d/2$, we pick a generic pencil of
hyperplane sections $H_t,\ t\in\Pee$ with smooth base-locus
$H_0\cap H_\infty\subset X$. Blowing this up gives another smooth
projective $d$-fold $\XX\stackrel{\pi\,}{\to}\Pee$ that fibres over
$\Pee$ with generically smooth fibres $H_t$. Fix any $t\in\Pee$
such that $H_t$ is smooth.

By the induction assumption there exist integers $N^t_1,\,N^t_2$ such that
$N^t_1A|^{\ }_{H_t}+N^t_2\omega^k$ is represented by an effective algebraic
cycle $Z_t\subset H_t$. Since, as Chad Schoen pointed out to me, the number
of such smooth fibres $H_t$ is uncountable, but the number of pairs
$(N_1,N_2)\in\Z\times\Z$ is countable, there exist an infinite number of
$t\in\Pee$ for which $(N^t_1,N^t_2)$ are all equal to some fixed $(N_1,N_2)$.

Consider now the relative Hilbert Scheme
\,Hilb\,$(\XX/\Pee)\to\Pee$ of cycles in the fibres $H_t$ of
cohomology class $N_1A+N_2\omega^k$. This is proper over $\Pee$
and surjects onto infinitely many points of $\Pee$, so surjects
onto $\Pee$. It is also projective, since $\XX$ is, so has a
degree $r$ multisection for some $r$. Pulling back the universal
subscheme via this multisection gives an effective algebraic cycle
$Z$ whose cohomology class on each fibre $H_t$ is
$r(N_1A+N_2\omega^k)$.

But by Lefschetz, the composition
$H^{2k}(X)\stackrel{\pi^*}{\longrightarrow} H^{2k}(\XX)\to
H^{2k}(H_t)$ is an isomorphism, so that the class of $\pi_*Z$ on
$X$ is $rN_1A+rN_2\omega^k$.
\end{Proof}

\section{If}

In this section we show that if the question in Theorem \ref{1} can be
answered positively then the statement in Proposition \ref{2} (and so the
Hodge conjecture) is true.

So we start with a rational class $A\in H^{n,n}(X^{2n})$ and assume there
is a nodal hypersurface $D\subset X$ and a homology class $B\in H_{2n}(D)$
whose image in $X$ is the dual of $A$.

Let $\{p_i\}\subset D$ denote the ODPs of $D$; blow these up inside $X$
to give $\pi:\,\Xx\to X$ with exceptional set $E=\cup_iE_i$ a collection
of $\mathbb P^{2n-1}$s.
Denote by $\D\stackrel{\pi\,}{\to}D$ the (smooth) proper transform of $D$
with exceptional set a collection of $(2n-2)$-dimensional quadrics
$Q_i\subset E_i$. We will also use $E$ and $E_i$ to denote the
$H^2$-classes of their corresponding line bundles, so that, on restriction
to $\D$, for instance, $E_i=Q_i$, etc.

\begin{Lemma}
$B$ lifts (non-canonically) to $\B\in H_{2n}(\D)$ such that $\pi_*\B=B$.
\end{Lemma}

\begin{Proof}
Choose a small contractible neighbourhood $U_i$ of each node $p_i\in D$,
with boundary $\partial U_i$ the link of the ODP. In $\D$, $\partial U_i$
bounds a tubular neighbourhood $N_i$ of the exceptional quadric $Q_i$.
$B$ defines a class in
$$
H_{2n}(D)\to H_{2n}(D,\cup_iU_i)\to H_{2n-1}(\cup_i\partial U_i)=
\bigoplus_i H_{2n-1}(\partial N_i),
$$
where the last group is surjected onto by $\bigoplus_i H_{2n}
(\cup_i N_i,\partial
N_i)$ since $H_{2n-1}(N_i)=H_{2n-1}(Q_i)=0$. This gives a relative homology
class to add to $B$ over the $Q_i$ to give $\B$ via Mayer-Vietoris.
\end{Proof}

It will now be sufficient to show that since $A$ is of pure Hodge
type, so is $\B$, since then we are reduced to finding a rational
combination of algebraic cycles to represent $\B$, which is the
Hodge conjecture in one dimension down (on $\D$), were our usual
induction argument applies. The pushdown to $D$ and
pushforward to $X$ of any such cycle in $\D$ would then be our
required cycle.

\begin{Lemma} $\B\in H_{2n}(\D)\cong H^{2n-2}(\D)$ is of Hodge type $(n-1,n-1)$.
\end{Lemma}

\begin{Proof}
Let $\iota:\,\D\to\Xx$ denote the inclusion. Then $\pi_*\iota_*\B$ is the
Poincar\'e dual of $A$, so that $\iota_*\B=\pi^*A\ $mod$\ \ker(\pi_*)$.
Both $A$ and the kernel of $\pi_*$ (consisting of exceptional classes on $\Xx$)
are pure of Hodge type $(n,n)$, so $\iota_*\B$ is too. Restricting its
cohomology
class back to $\D$ therefore shows that $\B\cup[\D]$ is of Hodge
type $(n,n)$ on $\D$, so it is sufficient to show that cupping
with $[\D]=\pi^*\omega-2[E]$ is an injection $H^{2n-2}(\D)\to
H^{2n}(\D)$, since it clearly preserves Hodge type.

Lemma 1.1 of \cite{Sch} gives the exact sequence
\begin{equation} \label{d}
0\to H^{2n-2}(X)\ \oplus\ \,\bigoplus_i\C\,.\,E_i^{n-1}\to H^{2n-2}(\D)
\to\langle A_i-B_i\rangle^*\to0, \vspace{-3mm}
\end{equation}
where the $A_i$ and $B_i$ are the standard $\mathbb P^{n-1}$
planes in the quadrics $Q_i$ (so that $E_i^n=A_i+B_i$ as
cohomology classes), and $\langle\ \rangle$ denotes their span in
$H^{2n}(\D;\C)$. Dually,
\begin{equation} \label{dd}
0\to\langle A_i-B_i\rangle\to H^{2n}(\D)\stackrel{\!\!\pi_*\iota_*}
{\longrightarrow}H^{2n+2}(X)\ \oplus\ \,\bigoplus_i\C\,.\,E_i^n\to0. \vspace{-2mm}
\end{equation}

The $\oplus_i\,\C\,.\,E_i^n$ part of (\ref{dd}) can be split by the obvious
restriction map $E_i^n\in H^{2n}(\Xx)\to H^{2n}(\D)$. To deal with the
remaining $H^{2n+2}(X)$ component, consider the composition $\alpha$,
$$
\alpha\,:\ H^{2n+2}(X)\cong H^{2n-2}(X)\to H^{2n}(X)\to H^{2n}(\D),
$$
of the Hard Lefschetz theorem followed by $\cup\,\omega$ then restriction to
$D$ and pullback to $\D$. This splits the sequence (\ref{dd}),
as it is easy to see that composing it with pushing forward
and back down to $X$ gives the identity. That is,
\begin{equation} \label{sp1}
H^{2n-2}(\D)\ \cong\ H^{2n-2}(X)\ \oplus\ \langle A_i-B_i\rangle^*\ \oplus\ \,
\bigoplus_i\C\,.\,E_i^{n-1}, \vspace{-3mm}
\end{equation}
and
\begin{equation} \label{sp2}
H^{2n}(\D)\ \cong\ H^{2n+2}(X)\ \oplus\ \langle A_i-B_i\rangle\ \oplus\ \,
\bigoplus_i\C\,.\,E_i^n. \vspace{-1mm}
\end{equation}
These splittings are compatible, in the sense that the corresponding terms
are duals of each other (via the obvious pairings on $X$ and $\D$ respectively)
and annihilate the other terms.

Since the dimensions match up, it is enough to show that the map $\theta_1=
\cup\,(\omega-2E):\,H^{2n-2}(\D)\to H^{2n}(\D)$ is onto. This map takes the
first summand of the splitting (\ref{sp1}) onto the first summand of (\ref{sp2})
(via the Lefschetz isomorphism $\cup\,\omega^2:\,H^{2n-2}(X)\to H^{2n+2}(X)$),
and the third summand onto the third summand (via the isomorphism
$\cup\,(-2E):\,\C\,.\,E_i^{n-1}\to\C\,.\,E_i^n$). The same is true of the
map $\theta_N=\cup\,(N\omega-2E)$, for any $N>0$. 

Let $\pi_2$ denote the projection to the second summand $\langle A_i-B_i\rangle$
of (\ref{sp2}); we are now left with showing that $\pi_2\theta_1$ is onto.
Since for $N\gg0$, $\OO(NH-2E)$ is an ample class on $\Xx$ (and so also on $\D$),
$\theta_N$ is onto. Therefore $\pi_2\theta_N$ is also onto, and it would be sufficient
to show that $\pi_2\theta_1=\pi_2\theta_N$, i.e. that $\pi_2\circ(\cup\,\omega)=0$.

But this follows from the fact that $\cup\,\omega:\,H^{2n-2}(\D)\to H^{2n}(\D)$
has image entirely in the first summand $H^{2n+2}(X)$, since it factors through
$$
H^{2n-2}(\D)\stackrel{\!\!\pi_*\iota_*}{\longrightarrow}H^{2n}(X)
\stackrel{\!\cup\omega}{\longrightarrow}H^{2n+2}(X)
\stackrel{\alpha\,}{\longrightarrow}H^{2n}(\D),
$$
as the generic hyperplane $H$ dual to $\omega$ misses the ODPs.
\end{Proof}

\begin{Remark}
Burt Totaro has pointed out to me that the result of this section holds
much more generally by using the machinery of Deligne's
mixed Hodge theory: any Hodge class on $X$ in the image of the pushfoward
on homology from a subvariety $D\subset X$ is the image of a Hodge class
on any resolution $\D$ of $D$.
\end{Remark}

\section{Only if}

To prove the only if part of Theorem \ref{1}, assume that the given rational
class $A\in H^{n,n}(X^{2n})$ is, up to taking rational multiples and adding
$n$-fold intersections of hyperplane sections of $X$, representable as an
effective algebraic cycle $Z$. Then by a theorem of Kleiman \cite{Kl} pointed
out to me by Daniel Huybrechts, $Z$ may be assumed smooth:

\begin{Theorem} \emph{\cite{Kl}} \label{kl}
Let $Z^k\subset X^d$ be an effective cycle of dimension
$k\le{d+1\over2}$ in a smooth projective $d$-fold $X$.
Then $(d-k-1)![Z]+NH^{\cap(d-k)}$ is algebraically equivalent to a
smooth cycle for $N\gg0$.
\end{Theorem}

Kleiman's result is proved as follows. Since $X$ is smooth and projective,
the structure sheaf $\OO_Z$ of $Z$ has
a finite locally free resolution of a standard form in which all of the
sheaves except for the last are of the form $E_i=
\OO(-n_iH)^{\oplus N_i}$ (inductively, set $E_i:=H^0(E_{i-1}(n_iH))\otimes
\OO(-n_iH)$ for sufficiently positive $n_i$). The $d$\,th sheaf $K$, the
(locally free) kernel of the resolution, therefore has Chern classes which
can be written
in terms of sums and products of the Chern classes of $\OO_Z$ and $[H]$.
In fact, up to multiples of $H^{\cap(d-k)}$, $c_{d-k}(K(NH))=(d-k-1)![Z]$
as cycles, so we may replace $Z$ by the zero section of the wedge of $d-k$
generic sections of some sufficiently positive twist $K(NH)$ of $K$.
Kleiman shows that for $N\gg0$ these sections may be taken to avoid certain
low dimensional Schubert cells (in the Grassmannian of $H^0(K(NH))$) which
parametrise degenerate $(d-k)$-tuples of sections; in this way he shows
the zero section can be taken to be smooth.

Now that $Z$ is smooth, we can show that it lies in a nodal hypersurface.

\begin{Theorem} \label{odp}
Let $Z^n\subset X^{2n}$ be a \emph{smooth} subvariety. For $N\gg0$, $Z$ is
contained in hypersurfaces $D$ in $|NH|$ whose only singularities are ODPs
on $Z$.
\end{Theorem}

\begin{Proof}
The base locus of the linear system $H^0(\I_Z(N))$ is $Z$ for $N\gg0$, so
the generic element is smooth away from $Z$. Restriction to $Z$ gives a map to
$H^0(\nu_Z^*(N))$ (where $\nu^*_Z=\I_Z/\I_Z^2$ is the conormal bundle to
$Z\subset X$), which is \emph{onto} for $N$ sufficiently large.
Sections of $\I_Z(N)$ vanish on $Z$, and this map takes their
derivative on $Z$. So the hyperplane is smooth along $Z$ where this
section of $\nu_Z^*(N)$ is non-zero, and has ODPs on $Z$ precisely where
it has simple transverse zeros. (Recall that a holomorphic
function $f$ on a smooth variety cuts out an analytic ordinary double point
at $x\in f^{-1}(0)$ if and only if $df$ has a simple zero at $x$.)

Since the rank of $\nu_Z$ is the same as the dimension $n$ of $Z$, for
$N\gg0$ the generic section of $\nu_Z^*(N)$ has a finite number
$c_n(\nu_Z^*(N))$ of simple zeros on $Z$, so the generic hyperplane
containing $Z$ has only ODPs on $Z$.
\end{Proof}

Thus, in particular, up to rational multiples and intersections of hyperplanes,
$A\in H^{2n}(X)$ is Poincar\'e dual to the pushforward to $X$ of a homology class
$[Z]\in H_{2n}(D)$ on $D$, as required.

\section{Infinitesimal Hodge conjecture}

This paper is the result of a failed attempt at proving Grothendieck's
variational Hodge conjecture; that if the class of an algebraic cycle $Z$
remains $(p,p)$ under an algebraic deformation of $X$, then its class
is also represented by a rational combination of algebraic cycles in the
deformation. We now explain why the method fails, as the obstruction is
interesting, and demonstrates that unfortunately, infinitesimally at least,
finding hyperplanes with ODPs is as hard as finding middle dimensional cycles.

The same induction as in Section 2 can be used for this
conjecture, though one cannot use the hard Lefschetz theorem to
pass from $k>d/2$ to $k<d/2$ this time (as we do not know that
there is an algebraic cycle $Z'$ whose intersection with $(2k-d)$
generic hyperplane sections is $Z$). Instead for such small cycles
one can first smooth using Theorem \ref{kl} and then use an easier
form of Theorem \ref{odp} (due to Altman and Kleiman \cite{AK}) to
include $Z$ in a \emph{smooth} hypersurface $D$. Some Lefschetz
theory then shows that as $X$ deforms (and $D$ with it, since the
deformation is assumed algebraic) that $Z\subset D$ remains of
type $(k-1,k-1)$ if $Z\subset X$ remains $(k,k)$. Thus we are
reduced to proving the variational Hodge conjecture for $Z\subset
D$ in one less dimension, and the induction procedure again leads
quickly to trying to prove the conjecture for middle dimensional
\emph{smooth} cycles $Z^n$ in \emph{even} dimensional varieties
$X^{2n}$.

So we include $Z$ inside a nodal hypersurface $D$, and try to show
that $D$ deforms with $X$ (preserving its ODPs), and $Z$ deforms
inside $D$. To show the
latter it would be sufficient to show that the proper transform
$\Zz$ of $Z$ in the blowup $\D$ of $D$ deforms; this can be shown
using our usual induction once we have shown that $\Zz\subset\D$
remains of type $(n-1,n-1)$ if $Z\subset X$ remains $(n,n)$ under
the deformation of $X$. With some work this can shown by Lefschetz methods
and results in \cite{Sch}, so we are reduced to showing that $D$ and its
ODPs do indeed deform with $X$.

The obstruction to deforming $D\in|NH|$ with $X$, with ODPs at the
points $p_i\in X$ (once we have chosen a fixed deformation of the
points $p_i$ with $X$), lies in $H^1(\I^2_{p_i}(NH))$. In fact we
now show that, allowing the $p_i$ to vary, instead of fixing their
deformation with $X$, the obstruction to deforming $D$ with the
ODPs preserved at any points lies in $H^1(\I_{p_i}(NH))$. The
latter group contains the obstruction to deforming $D$ while still
insisting it goes through the fixed $p_i$; if this obstruction
vanishes we get a section of $H^0(\I_{p_i}(NH))$ over the deformed
space, fitting into the exact sequence
$$
0\to H^0(\I^2_{p_i}(NH))\to H^0(\I_{p_i}(NH))\to\bigoplus_i
T^*_{p_i}X(NH).
$$
Thus the obstruction to lifting this to a hypersurface with ODPs
at the $p_i$ lies in the value of the derivative $ds|_{p_i}\in
T^*_{p_i}X(NH)$ of the given section $s$ of $\I_{p_i}(NH)$.
Deforming the $p_i$ by a vector $v_i\in T_{p_i}X$ changes this
derivative, and the corresponding differentiation map
$$
T_{p_i}X\otimes\big(s\in H^0(\I_{p_i}(NH))\big)\big|_{p_i}\to\ T^*_{p_i}X(NH)
$$
is a surjection, since by the definition of an ODP the Hessian of $s$ gives
an isomorphism $T_{p_i}X\to T^*_{p_i}X$ at each $p_i$.

Therefore, allowing for deformations of the $p_i$ we find that the obstruction
to deforming $D$ with ODPs at the $p_i$ lies in $H^1(\I_{p_i}(NH))$.

But now we find that, even as $N\to\infty$ (and so the number of $p_i$
tends to infinity too), this obstruction space does not disappear, in fact
it always contains the obstruction space $H^1(\nu_Z)$ to deforming $Z$!
The derivative $ds|_Z\in H^0(\nu_Z^*(NH))$ of the section $s$ of $\I_Z(NH)$
defining $D$ vanishes transversely at the $p_i$ as in Theorem \ref{odp},
giving the Koszul resolution on $Z$
\begin{equation} \label{kos}
0\to(\Lambda^n\nu_Z)((1-n)NH)\to\ldots
\to(\Lambda^2\nu_Z)(-NH)\to\nu_Z\to\I_{p_i}(NH)\to0.
\end{equation}
For $N$ sufficiently large this sequence and Kodaira vanishing give the
exact sequence
$$
0\to H^1(\nu_Z)\to H^1(\I_{p_i}(NH))\to H^n((\Lambda^n\nu_Z)((1-n)NH)),
$$
and so the inclusion claimed. Therefore deforming $D$ with
ODPs turns out to be as hard
as deforming $Z$. This is most easily seen in the apparently trivial $n=1$
case, where $Z\subset X$ is a curve in a surface, and $D$ is a nodal reducible
curve, the union of $Z$ and some other smooth curve. It is clear in this case
that to deform $D$ with these nodes deforms $Z$ too. \\

Finally we note, out of interest, that the extension class $e\in$\,Ext\,$^{n-1}
(\I_{p_i}(NH),$ $(\Lambda^n\nu_Z)((1-n)NH))$ of the sequence (\ref{kos}) lies
in
$$
\mathrm{Ext}^{\,n-1}_Z(\I_{p_i}(nD),\omega_Z\otimes\omega_X^*|^{\ }_Z)=
H^1_Z(\I_{p_i}(nD+\omega_X))^*,
$$
where $\omega$ denotes the canonical bundle. But for $N$ sufficiently large,
$H^1_Z(\I_{p_i}(nD+\omega_X))$ is the cokernel of $H^0_Z(\omega_X(nD)|_Z)\to
\bigoplus_i\omega_X(nD)|^{\ }_{p_i}$, which is the same as the cokernel
of $H^0_X(\omega_X(nD))\to\bigoplus_i\omega_X(nD)|^{\ }_{p_i}$, i.e.
$H^1_X(\I_{p_i}(nD+\omega_X))$.

But this last group is computed in \cite{Sch} to be the kernel of the pushforward
map $H^{2n}(\D)\to H^{2n+2}(\Xx)$. That is, the extension class of (\ref{kos})
lives in
$$
e\in{H^{2n-2}(\D)\over H^{2n-2}(\Xx)}\,,
$$
and is no doubt the fundamental class of $\Zz\subset\D$. (This group is
the $\langle A_i-B_i\rangle^*$ of (\ref{d}), i.e. the extra $2n$-homology
in $D$ that does not come from $X$ via the Lefschetz theorem for \emph{smooth}
hyperplanes. As $H^1_Z(\I_{p_i}(nD+\omega_X))^*$, it is easily computed
using (\ref{kos}) to be the one-dimensional
$H^n_Z(\Lambda^n\nu_Z\otimes\omega_X)^*\cong H^0(\OO_Z)$.)

It would be nice to try to use this in reverse; to at least identify (given
$A\in H^{2n}(X;\QQ)\cap H^{n,n}(X)$)
where the ODPs $p_i$ of the nodal hypersurface $D$ might lie in $X$ for
this relation to hold. This would be a first start in finding such a $D$
(and then $Z$), but seems unlikely to be possible as the above derivation
made such strong use of the cycle $Z$.

\begin{flushleft}
Department of Mathematics, Imperial College, London SW7 2AZ. UK. \\
\small Email: \tt rpwt@ic.ac.uk
\end{flushleft}


\begin{thebibliography}{Sch}

\bibitem[AK]{AK} A. Altman and S. Kleiman.
\emph{Bertini theorems for hypersurface sections containing a subscheme.}
Commun. Algebra \textbf{7}, 775--790 (1979). 

\bibitem[Cl]{Cl} H. Clemens.
\emph{Double solids.} Adv. in Math. \textbf{47}, 107--230 (1983).

\bibitem[GH]{GH} P. Griffiths and J. Harris.
\emph{Principles of algebraic geometry.} Wiley, New York, 1978.

\bibitem[Kl]{Kl} S. Kleiman.
\emph{Geometry on Grassmannians and applications to splitting bundles
and smoothing cycles.}
Publ. Math., Inst. Hautes \'Etud. Sci. \textbf{36}, 281--297 (1969).

\bibitem[Sch]{Sch} C. Schoen.
\emph{Algebraic cycles on certain desingularized nodal hypersurfaces.}
Math. Ann. \textbf{270}, 17--27 (1985).


\end{thebibliography}
\end{document}